\documentclass{amsart}
\usepackage{amsmath}
\usepackage{amssymb}
\usepackage{hyperref}
\hypersetup{
    linktocpage=true,
    colorlinks=true,
    linkcolor=blue,
    filecolor=magenta,      
    urlcolor=cyan,
    pdftitle={Extension of projection mappings},
    bookmarks=true,
   citecolor=red    }

\newtheorem{thm}{Theorem}[section]
\newtheorem{lem}[thm]{Lemma}
\newtheorem{prop}[thm]{Proposition}

\theoremstyle{definition}

\theoremstyle{remark}
\newtheorem{rem}[thm]{Remark}

\numberwithin{equation}{section}

\newcommand{\bh}{\mathcal{B}(H)} 

\newcommand{\tma}{S(\mathcal{A},\tau)}

\newcommand{\A}{\mathcal{A}}
\newcommand{\B}{\mathcal{B}}
\newcommand{\M}{\mathcal{M}}
\newcommand{\N}{\mathcal{N}}

\newcommand{\norm}[1]{\bigl\Vert{#1} \bigr\Vert}

\newcommand{\limm}[1]{\underset{#1}{\lim}}

\newcommand{\ra}{\rightarrow}

\newcommand{\supu}[1]{\underset{#1}{\sup}\,}

\newcommand{\id}{\mathbf{1}} 

\newcommand{\net}[1]{\{{#1}_{\lambda}\}_{\lambda \in \Lambda}}

\newcommand{\sotc}{\overset{SOT}{\rightarrow}}  

\newcommand{\seq}[1]{({#1}_n)_{n=1} ^{\infty}}

\newcommand{\summ}[2]{\underset{#1}{\overset{#2}{\sum}}}

\newcommand{\ip}[2]{\left\langle {#1},{#2} \right\rangle}

\newcommand{\rax}[1]{\overset{#1}{\rightarrow}}

\newcommand{\raxx}[1]{\underset{#1}{\rightarrow}}

\newcommand{\limx}[1]{\underset{#1}{\lim}}
\newcommand{\pa}{\mathcal{P}(\mathcal{A})}    
\newcommand{\pb}{\mathcal{P}(\mathcal{B})}    
\newcommand{\paf}{\mathcal{P}(\mathcal{A})_f}    
\newcommand{\pbf}{\mathcal{P}(\mathcal{B})_f}    

\newcommand{\mtop}{\mathcal{T}_m}

\newcommand{\D}{\mathcal{D}}
\newcommand{\xmx}[2]{S({#2},{#1})} 
\newcommand{\nmb}{S(\mathcal{B},\nu)} 
\newcommand{\n}{\mathbb{N}}

\newcommand{\Dx}[1]{\mathcal{D}_{(#1)}}
\newcommand{\uparrowx}[1]{\uparrow_{#1}}
\newcommand{\sotcx}[1]{\overset{SOT}{\underset{#1}{\rightarrow}}}  
\newcommand{\G}{\mathcal{G}}
\newcommand{\jnxx}[2]{\underset{#1}{\overset{#2}{\vee}}}
\newcommand{\netp}[1]{\{{#1}_{\lambda'}\}_{\lambda' \in \Lambda'\subseteq \Lambda}}

\begin{document}

\title[Extension of projection mappings]{Extension of projection mappings}


\author{Pierre de Jager}
\address{Department of Mathematics, University of Cape Town, South Africa}
\curraddr{DST-NRF CoE in Math. and Stat. Sci\\ Unit for BMI\\ Internal Box 209, School of Comp., Stat., $\&$ Math. Sci.\\
NWU, PVT. BAG X6001, 2520 Potchefstroom\\ South Africa}
\email{28190459@nwu.ac.za}
\thanks{The first author would like to thank the NRF for funding towards this project in the form of scarce skills and grantholder-linked bursaries}

\author{Jurie Conradie}
\address{Department of Mathematics, University of Cape Town, South Africa}
\email{jurie.conradie@uct.ac.za}


\subjclass[2010]{Primary 46L10; Secondary 46L51}
\keywords{Extensions, Jordan homomorphisms, semi-finite von Neumann algebras}

\date{\today}




\maketitle

\begin{abstract}
We show that a map between projection lattices of semi-finite von Neumann algebras can be extended to a Jordan $*$-homomorphism between the von Neumann algebras if this map is defined in terms of the support projections of  images (under the linear map) of projections and the images of orthogonal projections have orthogonal support projections. This has numerous fundamental applications in the study of isometries and composition operators on quantum symmetric spaces and is of independent interest, since it provides a partial generalization of Dye's Theorem without the requirement that the initial von Neumann algebra be free of type $I_2$ summands.
\end{abstract}

\section{Introduction}

In \cite{key-Y81}, Yeadon showed that an isometry between quantum (non-commutative) $L^p$-spaces associated with semi-finite von Neumann algebras can essentially be characterized as a weighted non-commutative composition operator. A fundamental component of the proof of this result is using the support projections of the images (under the isometry) of projections of finite trace to define a mapping of projections that can be extended to a Jordan $*$-homomorphism. This extension procedure uses both the linearity of the isometry and aspects of the $L^p$-space structure. In order to characterize isometries in other settings it is desirable to  have a more general framework for these extension procedures. We will therefore be interested in describing conditions under which a mapping between projection lattices of von Neumann algebras (we will call such a map a \emph{projection mapping}) can be extended to a Jordan $*$-homomorphism between the von Neumann algebras. It will be shown in \cite{key-dJ1}, \cite{key-dJ2} and \cite{key-dJ3} that the extension procedures developed herein can be used in the characterizations of positive surjective isometries on quantum symmetric spaces, surjective isometries on Lorentz spaces and isometries on Orlicz spaces respectively. (These results form part of the first author's thesis \cite{key-dJ17}.) Furthermore, the procedures developed herein will also facilitate the characterization of quantum composition operators between symmetric spaces (\cite{key-dJ18c}), which requires the extension of a projection-preserving mapping from a certain subalgebra of a von Neumann algebra to the whole algebra.

It is important to note the work of Dye, who showed (\cite{key-Dye55}) that a projection ortho-isomorphism (bijective mapping between projection lattices of von Neumann algebras preserving orthogonality in both directions) is necessarily implemented by a Jordan $*$-isomorphism. Unfortunately, projection ortho-isomorphisms have significantly more structure than the projection mappings typically obtained in the characterization of isometries. It is also worth noting that in \cite{key-Bunce93} it is shown that any map between projection lattices of von Neumann algebras, where the domain von Neumann algebra has no type $I_2$ direct summand, is a restriction of a Jordan $*$-homomorphism if and only if it is additive on orthogonal projections. An important related topic is the Mackey-Gleason problem, which deals with determining conditions under which bounded vector measures can be extended to bounded linear maps (see \cite{key-Bunce92} and \cite{key-Bunce94}, for example, or the introduction to Hamhalter's article (\cite{key-Ham15}) for an overview and some historical background). Significantly and unsurprisingly solutions of the Mackey-Gleason problem are also confined to the setting where the initial von Neumann algebra is free of type $I_2$ summands. Since we do not wish to have this restriction and with our particular applications in mind, we will proceed with the special case where the projection mapping is given in terms of the support projections of images of projections under a linear mapping. We will also confine ourselves to semi-finite von Neumann algebras. Before describing in more detail the aims and structure of this paper, we first introduce the notation we will be using and provide some background information.

\section{Preliminaries}

Throughout this paper we will use $\A\subseteq B(H)$ and $\B\subseteq B(K)$ to denote semi-finite von Neumann algebras, equipped with faithful normal semi-finite traces $\tau$ and $\nu$ respectively. Let $\pa$ denote the lattice of projections in $\A$, and $\paf$ the sublattice consisting of projections with finite trace. The set of all finite linear combinations of mutually orthogonal projections in $\pa$ (alternatively $\paf$) will be denoted $\mathcal{G}$ (respectively $\mathcal{G}_f$). We will use $x_\alpha \ra x$, $x_\alpha \rax{SOT} x$ and $x_\alpha \rax{WOT} x$ to denote convergence with respect to the norm topology, strong operator topology (SOT) and weak operator topology (WOT) respectively. We will repeatedly use the facts that addition is jointly SOT-continuous, and multiplication is jointly SOT-continuous provided one of the variables is restricted to a bounded set. In these situations we will typically omit the details regarding the indexing sets and write $x_\alpha+y_\beta \sotc x+y$ (or $x_\alpha y_\beta \sotc xy$) if $x_\alpha \sotc x$ and $y_\beta \sotc y$. Further details regarding von Neumann algebras may be found in \cite{key-K1}.

The set of all trace-measurable operators affiliated with $\A$ and $\B$ will be denoted $\tma$ and $\nmb$, respectively. Equipped with the measure topology $\mtop$, these spaces are complete metrisable  topological $*$-algebras. If $\mathcal{H}$ is a collection of trace-measurable operators, then we will use $\mathcal{H}^{sa}$ to denote the self-adjoint elements in $\mathcal{H}$ and $\mathcal{H}^+$ to denote the positive ones. If $\{x_\lambda\}_{\lambda \in \Lambda}\cup\{x\}$ is a collection of self-adjoint trace-measurable operators, then we will write $x_\lambda \uparrow x$ if $\{x_\lambda\}_{\lambda \in \Lambda}$ is an increasing net and $x=\supu{\lambda}x_\lambda$. If $\mathcal{H}\subseteq \tma$ and $T:\mathcal{H} \ra \nmb$ is a linear map such that $T(x_\lambda)\uparrow T(x)$ whenever $\{x_\lambda\}_{\lambda \in \Lambda}\cup\{x\}\subseteq \mathcal{H}^{sa}$ is such that $x_\lambda \uparrow x$, then $T$ will be called \emph{normal} (on $\mathcal{H}$). The support and range projections of a trace-measurable operator $x$ will be denoted $s(x)$ and $r(x)$, respectively. For background and further details regarding trace-measurable operators the interested reader is directed to \cite{key-Dodds14}.

The following consequence of the spectral theorem is well-known and will play a significant role in our extension procedures.

\begin{prop} \label{RP1 19/11/15}
If $x\in \A^{sa}$, then there exits a sequence $\seq{x}\subseteq \G^{sa}$ such that $x_n \ra x$; $s(x_n)\leq s(x)$ for all $n\in \mathbb{N}^+$ and the $x_n$'s commute with each other and with $x$. If, in addition $\tau(r(x))<\infty$, then $x_n \in \G_f^{sa}$ for every $n\in \mathbb{N}^+$. If $x\geq 0$, then $x_n\geq 0$ for every $n\in \mathbb{N}$. 
\end{prop}

 A linear mapping $\Phi :\M \ra \N$ between $C^*$-algebras is called a \emph{Jordan homomorphism} if $\Phi(yx+xy)=\Phi(y)\Phi(x)+\Phi(x)\Phi(y)$ for all $x,y \in \M$. If, in addition, $\Phi(x^*)=\Phi(x)^*$ for all $x\in \M$, then $\Phi$ is called a \emph{Jordan $*$-homomorphism}. 

\begin{rem} \label{alt def Jordan}
Some authors define a Jordan $*$-homomorphism to be $*$-preserving linear map such that $\Phi(x^2)=\Phi(x)^2$ for all $x\in \M^{sa}$. Noting that $(x+y)^2=x^2+xy+yx+y^2$ for all $x,y\in \M$ and recalling that any element $x\in \M$ can be written in the form $x=x_1+ix_2$, where $x_1,x_2 \in \M^{sa}$, it is easily checked that these definitions are equivalent.
\end{rem}

We will often be interested in maps with similar properties to Jordan $*$-homomor-phisms, but defined on $\mathcal{F}(\tau)=\{x\in \A: \tau(s(x))<\infty\}$, which is a (not necessarily closed) subalgebra of the von Neumann algebra $\A$. The following details properties of such maps.

\begin{prop} \label{P18 16/10/14} 
If $\Phi: \mathcal{F}(\tau) \ra \B$ is a positive linear map such that $\Phi(x^2)=\Phi(x)^2$ for all $x\in \mathcal{F}(\tau)^{sa}$, then, for all $x,y\in \mathcal{F}(\tau)$,
\begin{enumerate}
\item $\Phi(x^*)=\Phi(x)^*$;
\item $\Phi(xy+yx)=\Phi(x)\Phi(y)+\Phi(y)\Phi(x)$;
\item $\Phi(x^n)=\Phi(x)^n$ for all $n\in \mathbb{N}^+$;
\item $\Phi(xyx)=\Phi(x)\Phi(y)\Phi(x)$;
\item $\Phi(p)$ is a projection for every $p\in \paf$; 
\item $\Phi(p)\Phi(q)=0$ whenever $p,q\in \paf$ with $pq=0$;
\item $\Phi(x)\Phi(y)=\Phi(y)\Phi(x)$ and $\Phi(xy)=\Phi(x)\Phi(y)$ whenever $xy=yx$;
\item $\norm{\Phi(x)}_\B \leq \norm{x}_\A$ whenever $x\in \mathcal{F}(\tau)^{sa}$.
\end{enumerate}
\begin{proof}
It is easily checked that (1) follows from the fact that $\Phi$ is linear and positive. A similar argument to the one used in Remark \ref{alt def Jordan} shows that (2) holds. Properties (3) to (7) are algebraic in nature and can therefore be shown using the same techniques as in the proofs of the corresponding properties in \cite[Exercises 10.5.21 and 10.5.22]{key-K4}. To prove (8), let $x\in  \mathcal{F}(\tau)^{sa}$. It follows by \cite[Proposition 4.2.3]{key-K1} that $-\norm{x}\id\leq x\leq \norm{x}\id$. Since $x$ is self-adjoint, $s(x)xs(x)=x$, and therefore \[-\norm{x}s(x) \leq x \leq \norm{x}s(x),\] by \cite[Proposition 1(iii)]{key-Dodds14}. Since $x\in \mathcal{F}(\tau)^{sa}$, $\tau(s(x))<\infty$ and therefore $s(x)\in \mathcal{F}(\tau)$. Since $\Phi$ is positive and linear we therefore have \[-\norm{x}\Phi(s(x)) \leq \Phi(x) \leq \norm{x}\Phi(s(x)).\] By \cite[Proposition 4.2.8]{key-K1}, this implies that $\norm{\Phi(x)}\leq \norm{x}\norm{\Phi(s(x))}$. This completes the proof, since $\Phi(s(x))$ is a projection (using (5)) and therefore $\norm{\Phi(s(x))}=1$.
\end{proof}
\end{prop}

Suppose $\Phi$ is a map from $\paf$ into $\pb$.  Ideally, one would like to know the conditions on $\Phi$  which would ensure that it is uniquely extensible to a Jordan $*$-homomorphism from $\A$ into $\B$ (or preferably, in certain instances, from $\A$ onto $\B$). There appears to be several difficulties. The first is in ensuring that the extension to linear combinations of projections is linear (and not just linear on commuting elements). This can be overcome if there exists a linear map $U:\mathcal{F}(\tau) \ra \xmx{\nu}{\B}$ with the property that $\Phi(p)=s(U(p))$ for all $p\in \paf$. This condition is restrictive, but is sufficiently general for the extension procedure based on it to be useful in the characterization of isometries between various non-commutative spaces.  The second difficulty is in extending $\Phi$ from $\mathcal{F}(\tau)$ to $\A$, which appears to require normality of the map $\Phi$ on $\mathcal{F}(\tau)$ to ensure that the extension to $\A$ is linear. Motivated by this and also the desire to determine conditions for the extensibility of a projection-preserving mapping on $\mathcal{F}(\tau)$, we have divided the extension process into several steps. In section \ref{S3} we will describe conditions on a projection mapping on $\paf$ which ensures that it can be extended to a map on $\mathcal{F}(\tau)$ with Jordan $*$-homomorphism-like properties. In section \ref{S4} such a map will be extended to a Jordan $*$-homomorphism on the whole von Neumann algebra; this can be done if the map is also normal. It is worth noting that in the study of isometries it is often the case that the interaction of the properties of the isometry and the spaces under consideration can be used to show the normality of the extension obtained in section \ref{S3} and hence ensure its further extension. We will, however, also consider conditions on a projection mapping which ensures its extension to a Jordan $*$-homomorphism (see section \ref{S5}). In the final section we consider sufficient conditions to ensure the surjectivity of a Jordan $*$-homomorphism.

\section{Extending projection mappings to $\mathcal{F}(\tau)$}\label{extension2}\label{S3}

In this section we show that under certain conditions a projection mapping $\Phi$ defined on $\paf$ can be extended to a positive square-preserving linear contraction on $\mathcal{F}(\tau)$. The proof is fairly lengthy and so we have divided it into a number of lemmas. We start by showing that a projection mapping on an arbitrary lattice of projections can be extended to the set of  finite real linear combinations of those projections. This will be used for the extension of a map defined on the set of all projections of finite trace in the present section and also for the extension of a map defined on the set of all projections in the sequel. Surprisingly and unfortunately this result only guarantees the linearity of this extension on commuting elements. We will therefore prove a sequence of lemmas that will enable us to show the linearity of this extension on all elements in its domain, in the case where $\Phi(p)=s(U(p))$, for all $p\in \paf$, for some linear map $U:\mathcal{F}(\tau)\ra \nmb$. This will enable us to further extend this map to a complex linear map on all of $\mathcal{F}(\tau)$. We will then show that this extension has the desired properties.

\begin{lem} \label{nL1 04/12/15} 
Suppose $\mathcal{E}\subseteq \pa$ is a lattice of projections containing the zero projection and let $\mathcal{H}$ denote the set of all finite linear combinations of mutually orthogonal projections from $\mathcal{E}$. If $\Phi:\mathcal{E}\ra \pb$ is a map such that $\Phi(p+q)=\Phi(p)+\Phi(q)$, whenever $p,q\in \mathcal{E}$ with $pq=0$, then $\Phi(0)=0$, $\Phi(p)\Phi(q)=0$ whenever $p,q\in \mathcal{E}$ with $pq=0$, and $\Phi(p)\geq \Phi(q)$ whenever $p,q\in \mathcal{E}$ with $p\geq q$. Furthermore, $\Phi$ can be extended to a  map (still denoted $\Phi$) from $\mathcal{H}^{sa}$ into $\B$, which is real linear on commuting elements, and has the properties $\Phi(x^2)=\Phi(x)^2$ and $\norm{\Phi(x)}_\B\leq\norm{x}_\A$ for all $x\in \mathcal{H}^{sa}$.
\begin{proof}
It is easily checked that the additivity of $\Phi$ yields $\Phi(0)=0$, $\Phi(p)\Phi(q)=0$, whenever $p,q\in \mathcal{E}$ with $pq=0$, and that $\Phi(p)\geq \Phi(q)$ whenever $p,q\in \mathcal{E}$ with $p\geq q$. For $0\neq x=\summ{i=1}{k}\alpha_i p_i \in \mathcal{H}^{sa}$,  let $\Phi(x):=\summ{i=1}{k}\alpha_i \Phi(p_i)$. We show that $\Phi$ is well-defined. Suppose $x=\summ{j=1}{m}\beta_j q_j$ is another representation of $x$.  Note that $\summ{i=1}{k}p_i=\jnxx{i=1}{k}p_i\in \mathcal{E}$ and $\summ{j=1}{m}q_j =\jnxx{j=1}{m}q_j\in \mathcal{E}$, since these are sums of mutually orthogonal projections. Furthermore, $\summ{i=1}{k}p_i=\summ{j=1}{m}q_j$, since we can assume, without loss of generality, that $\alpha_i\neq 0$ for all $i$ and $\beta_j\neq 0$ for all $j$. These sums are therefore representations of the support projection of $x$. This implies that $\summ{i=1}{k}\Phi(p_i)=\Phi(\summ{i=1}{k}p_i)=\Phi(\summ{j=1}{m}q_j)=\summ{j=1}{m}\Phi(q_j)$. Since, for $i\neq j$, $p_j p_i=0$  implies that $\Phi(p_j)\Phi(p_i)=0$ and therefore $\summ{j=1}{m}\Phi(q_j)\in \pb$ and $\Phi(p_i)\leq \summ{j=1}{m}\Phi(q_j)$, we have that, for all $1\leq i\leq n$,
\begin{eqnarray}
\Phi(p_i)&=&(\summ{j=1}{m}\Phi(q_j))\Phi(p_i)=\summ{j=1}{m}(\Phi(q_j)\Phi(p_i)). \label{e13.1.0 C1}
\end{eqnarray}
We can similarly show that for every $1\leq j \leq m$,
\begin{eqnarray}
\Phi(q_j)=\Phi(q_j)(\summ{i=1}{k}\Phi(p_i))=\summ{i=1}{k}\left(\Phi(q_j)\Phi(p_i)\right). \label{e13.1.1 C1}
\end{eqnarray}
Furthermore, if $1\leq i'\leq k$ and $1\leq j' \leq m$, it follows that 
\begin{eqnarray*}
\alpha_{i'} q_{j'} p_{i'}=q_{j'} (\summ{i=1}{k}\alpha_i p_i)p_{i'}=q_{j'}(\summ{j=1}{m}\beta_j q_j)p_{i'} = \beta_{j'} q_{j'} p_{i'} 
\end{eqnarray*}
and hence that $\alpha_{i'}=\beta_{j'}$ or $q_{j'} p_{i'}=0$. Since $\Phi$ maps orthogonal projections onto orthogonal projections, we have that for every $i,j$, either $\alpha_i=\beta_j$ or $\Phi(q_j)\Phi( p_i) = 0$. Therefore, using (\ref{e13.1.0 C1}) and (\ref{e13.1.1 C1}), we obtain
\begin{eqnarray*}
\summ{i=1}{k}\alpha_i \Phi(p_i)&=&\summ{i=1}{k}\alpha_i \summ{j=1}{m}\Phi(q_j)\Phi(p_i) =\summ{i=1}{k}\summ{j=1}{m}\alpha_i \Phi(q_j)\Phi(p_i)  =\summ{i=1}{k}\summ{j=1}{m}\beta_j \Phi(q_j)\Phi(p_i)\\
&=&\summ{j=1}{m}\summ{i=1}{k} \beta_j \Phi(q_j)\Phi(p_i) =\summ{j=1}{m}\beta_j \summ{i=1}{k} \Phi(q_j)\Phi(p_i) =\summ{j=1}{m}\beta_j\Phi(q_j).
\end{eqnarray*}

It follows that $\Phi$ is well-defined on $\mathcal{H}^{sa}$. Next, we show that $\Phi$ is real linear on commuting operators in $\mathcal{H}^{sa}$. It clearly suffices to show that $\Phi$ is additive. Suppose $x=\summ{i=1}{k}\alpha_ip_i$ and $y=\summ{j=1}{m}\beta_jq_j$ are commuting operators in $\mathcal{H}^{sa}$. Put $p=\summ{i=1}{k}p_i$ and $q=\summ{j=1}{m}q_j$. Define $r_{i,j}:=p_iq_j$, $s_i:=p_i-p_iq$ and $t_j:=q_j-q_jp$. Note that since we can assume, without loss of generality, that $\alpha_i \neq \alpha_j$ and $\beta_i\neq \beta_j$ if $i\neq j$, the $p_i$'s and $q_i$'s are spectral projections of $x$ and $y$, respectively. Since $xy=yx$, we have that $p_iy=yp_i$ for all $i$, by \cite[Theorem 1.5.12]{key-dP} and \cite[Corollary 15]{key-Terp1}. It therefore follows by the same corollary that $p_iq_j=q_jp_i$ for all $i,j$ and so $r_{i,j},s_i$ and $t_j$ are projections for all $i,j$. Furthermore, these are mutually orthogonal, since $(p_i)_{i=1}^k$ and $(q_j)_{j=1}^m$ are families of mutually orthogonal projections. It is easily verified that
\begin{eqnarray*}
x+ y &=& \summ{i=1}{k}\summ{j=1}{m}(\alpha_i+ \beta_j)r_{i,j}+\summ{i=1}{k} \alpha_i s_i+\summ{j=1}{m} \beta_j t_j \in \mathcal{H}^{sa}. 
\end{eqnarray*}
Furthermore, $p_i=s_i+\summ{j=1}{m}r_{i,j}$ implies that $\Phi(p_i)=\Phi(s_i)+\summ{j=1}{m}\Phi(r_{i,j})$, since $\Phi$ is additive on orthogonal projections. Similarly, $\Phi(q_j)=\Phi(t_j)+\summ{i=1}{k}\Phi(r_{i,j})$. A straightforward calculation using these facts and grouping appropriate terms, yields $ \Phi(x)+\Phi(y)=\Phi( x+ y)$.

We show that $\Phi$ is a square-preserving contraction. If $x=\summ{i=1}{n}\alpha_i p_i \in \mathcal{H}^{sa}$, then $x^2=(\summ{i=1}{n}\alpha_ip_i)^2 =\summ{i=1}{n}\alpha_i^2p_i\in \mathcal{H}^{sa}$. Therefore, using the fact that $\Phi(p_i)\Phi(p_j)=0$ for $i\neq j$, we obtain
\[\Phi(x^2)=\summ{i=1}{n}\alpha_i^2 \Phi(p_i) =(\summ{i=1}{n}\alpha_i \Phi(p_i))^2 =\Phi(x)^2\] 
and
\[\norm{x}_{\A}=\norm{\summ{i=1}{n}\alpha_i p_i}_{\A} =\text{max}\{|\alpha_i|:i=1,2,...,n\} \geq \norm{\summ{i=1}{n}\alpha_i\Phi(p_i)}_{\B} =\norm{\Phi(x)}_{\B}.\]
\end{proof}
\end{lem}

\begin{rem}\label{R1 01/07/18}
If, in addition to the assumptions of Lemma \ref{nL1 04/12/15}, we have that $p=0$ whenever $\Phi(p)=0$, then $\Phi$ will be isometric on $\mathcal{H}^{sa}$. 
\end{rem}

For the remainder of this section, we will suppose that $\Phi:\paf \ra \pb$ is a map such that 
\begin{itemize}
\item $\Phi(p+q)=\Phi(p)+\Phi(q)$ whenever $p,q\in \paf$ with $pq=0$; 
\item $\Phi(p)=s(U(p))$ for all $p\in \paf$, for some linear map $U: \mathcal{F}(\tau) \ra\xmx{\nu}{\B}$ with the property that $U(x_n)\rax{\mtop} U(x)$ whenever $\seq{x}\cup \{x\}\subseteq \mathcal{F}(\tau)$ is such that $x_n \ra x$ and $s(x_n)\leq s(x)$ for all $n\in \mathbb{N}^+$. 
\end{itemize}

\begin{rem}
It is natural to ask if such projection mappings exist. The conditions listed above are motivated by the proof of Yeadon's Theorem (\cite[Theorem 2]{key-Y81}), where an isometry $U$ between non-commutative $L^p$-spaces is used to define a projection-mapping $\Phi$ by letting $\Phi(p)=s(U(p))$. It will be shown (\cite{key-dJ1}, \cite{key-dJ2} and \cite{key-dJ3}) that the projection mappings considered here arise naturally in the study of isometries between various types of non-commutative spaces.
\end{rem}

Since $\paf$ is a lattice and $\Phi$ satisfies the conditions of Lemma \ref{nL1 04/12/15}, $\Phi$ can be extended to a map $\Phi_1:\mathcal{G}_f \ra \B$, which is square-preserving on $\mathcal{G}_f^{sa}$, real linear on commuting elements and a contraction on $\mathcal{G}_f^{sa}$. It follows from Proposition \ref{RP1 19/11/15} that for any $x\in \mathcal{F}(\tau)^{sa}$  there exists a sequnce $(x_n)_{n=1}^\infty\subseteq \mathcal{G}_f^{sa}$ such that $x_n\ra x$. Since $\B$ is complete and $\Phi_1$ a contraction, which is real linear on commuting elements, letting $\Phi_2(x):=\limm{n\ra \infty} \Phi_1(x_n)$ yields a well-defined contraction $\Phi_2: \mathcal{F}(\tau)^{sa}\ra \B$. Unfortunately, Lemma \ref{nL1 04/12/15} only guarantees real linearity of $\Phi_1$ on commuting self-adjoint elements and hence also the real linearity of $\Phi_2$ is restricted to commuting self-adjoint elements. We will prove two lemmas (using the linear map $U$) that will enable us to show the linearity of $\Phi_2$ on all elements in its domain. To simplify notation, we will use $\Phi$ to denote both of the extensions $\Phi_1$ and $\Phi_2$. 

\begin{lem} \label{L15.1 C1} If $x\in \mathcal{F}(\tau)^{sa}$, then $U(x)=U(p)\Phi(x)$ for any $p\in \paf$ with $p\geq s(x)$.
\begin{proof} Suppose $x=q$ for some $q\in \paf$. If $p\in \paf$ with $p\geq s(q)=q$, then $p-q\in \paf$ and $s(U(p-q))= \Phi(p-q)$. Since $q(p-q)=0$ implies $\Phi(q)\Phi(p-q)=0$, it follows that $U(p-q)\Phi(q)=0$. Using this and the fact that $\Phi(q)=s(U(q))$, we obtain
\begin{eqnarray}
U(p)\Phi(x)=(U(p-q)+U(q))\Phi(q) =U(q)\Phi(q) =U(x).  \label{e15.2 C1}
\end{eqnarray}
Using the real linearity of $\Phi$ on commuting elements in $\mathcal{G}_f^{sa}$, we can extend (\ref{e15.2 C1}) to $x \in \mathcal{G}^{sa}$ and $p\in \paf$ with $p\geq s(x)$. For a general $x\in \mathcal{F}(\tau)^{sa}$ and $p\in \paf$ with $p\geq s(x)$, we can use Proposition \ref{RP1 19/11/15}, to find a sequence $\seq{x} \subseteq \mathcal{G}_f^{sa}$ such that $x_n\ra x$, $s(x_n)\leq s(x)\leq p$ for every $n\in \n^+$,  and the $x_n$'s commute with each other and with $x$. Since $\Phi$  is real linear on commuting elements and is a contraction, it follows that  $\Phi(x_n)\ra \Phi(x)$. Therefore $U(p)\Phi(x_n)\ra U(p)\Phi(x)$ and hence $U(p)\Phi(x_n)\rax{\mtop} U(p)\Phi(x)$, by \cite[Proposition 20]{key-Dodds14}. However, $U(p)\Phi(x_n)=U(x_n)$ for all $n\in \n^+$, by what has been shown already. It follows that $U(x_n)\rax{\mtop} U(p)\Phi(x)$. However, by the assumption on $U$ we also have that $U(x_n)\rax{\mtop} U(x)$. Since the measure topology is Hausdorff, we obtain $U(x)=U(p)\Phi(x)$. 
\end{proof}
\end{lem}

\begin{lem} \label{L15.2 C1}
If $x\in \mathcal{F}(\tau)^{sa}$, then $r(\Phi(x))\leq s(U(p))=\Phi(p)$ for any $p\in \paf$ with $p\geq s(x)$.
\begin{proof}
Suppose $x=q\in \paf$ and $p\in \paf$ with $p\geq s(x)= q$. This implies that $s(U(p))=\Phi(p)\geq \Phi(q)=r(\Phi(q))$, by Lemma \ref{nL1 04/12/15}. If $x=\summ{i=1}{n}\alpha_i p_i \in \mathcal{G}_f^{sa}$ and $p\geq s(x)=\summ{i=1}{n} p_i$, then 
\begin{eqnarray}
&&r(\Phi(x))=r(\summ{i=1}{n}\alpha_i \Phi(p_i)) =\summ{i=1}{n} \Phi(p_i) =\Phi(\summ{i=1}{n} p_i) \leq\Phi(p) =s(U(p)), \label{e15.2.5 C1}
\end{eqnarray}
since the $\Phi(p_i)$'s are orthogonal projections, $\Phi$ is additive on orthogonal projections and $p\geq \summ{i=1}{n} p_i$. Suppose $x\in \mathcal{F}(\tau)^{sa}$ and $p\in \paf$ with $p\geq s(x)$. By Proposition \ref{RP1 19/11/15}, there exists a sequence $\seq{x} \subseteq \mathcal{G}_f^{sa}$ such that $x_n \ra x$, $s(x_n)\leq s(x)\leq p$ for all $n\in \n^+$, and the $x_n$'s commute with each other and with $x$. It follows by (\ref{e15.2.5 C1}) that $r(\Phi(x_n))\leq s(U(p))$ for all $n\in \n^+$. Furthermore, $\Phi(x_n)\ra \Phi(x)$, since $\Phi$ is real linear on commuting elements and a contraction. For a trace-measurable operator $y$, let $\ker(y)$ denote the projection onto the kernel of $y$. Suppose $\eta$ is in the kernel of $(U(p))$. Note that for any $n \in \n^+$,
\begin{eqnarray}
&&\ker(U(p)) =[s(U(p))]^{\perp} \leq [r(\Phi(x_n))]^{\perp}=[r(\Phi(x_n))^*)]^{\perp}=\ker(\Phi(x_n)), \label{e15.2.6 C1}
\end{eqnarray}
using \cite[Exercise 2.8.45]{key-K3} and the facts that $r(\Phi(x_n))\leq s(U(p))$ (see (\ref{e15.2.5 C1})) and that $\Phi(x_n)$ is a real linear combination of projections. Furthermore, $\Phi(x_n)\ra \Phi(x)$ implies that $\Phi(x_n)\sotc \Phi(x)$ and so $\norm{\Phi(x)(\eta)}=\limx{n\ra \infty} \norm{\Phi(x_n)(\eta)} =0$, since $\eta$ is in the kernel of $(\Phi(x_n))$ for every $n \in \n^+$, using (\ref{e15.2.6 C1}). Furthermore, since $\B^{sa}$ is closed in $\B$ and $\Phi(\mathcal{G}_f^{sa})\subseteq \B^{sa}$, we have that $\Phi$ also maps into $\B^{sa}$. Using these facts and \cite[Exercise 2.8.45]{key-K3}, it follows that $\eta \in \ker(\Phi(x))(K)=[r(\Phi(x)^*)]^{\perp}(K)=[r(\Phi(x))]^{\perp}(K)$. Therefore $s(U(p))^{\perp}=\ker(s(U(p))) \leq r(\Phi(x))^{\perp}$ and hence $r(\Phi(x)) \leq s(U(p))$.
\end{proof}
\end{lem}

We are now in a position to show that $\Phi$ is real linear on $\mathcal{F}(\tau)^{sa}$. Suppose $x,y\in \mathcal{F}(\tau)^{sa}$ and $\alpha,\beta\in \mathbb{R}$. Let $p=s(x)\vee s(y)$. Then using Lemma \ref{L15.1 C1} and the linearity of $U$, we have 
\[U(p)[\Phi(\alpha x +\beta y)-\alpha \Phi(x) - \beta \Phi(y)]=U(\alpha x +\beta y)-\alpha U(x) - \beta U(y)  =0.\] 
Furthermore,
\[r([\Phi(\alpha x +\beta y)-\alpha \Phi(x) - \beta \Phi(y)])\leq r(\Phi(\alpha x +\beta y))\vee r(\Phi(x))\vee r(\Phi(y)) \leq s(U(p)), \]
by Lemma \ref{L15.2 C1}, since $s(x),s(y),s(\alpha x+\beta y)\leq p$.  It follows that $[\Phi(\alpha x +\beta y)-\alpha \Phi(x) - \beta \Phi(y)]=0$ and therefore $\Phi$ is real linear on $\mathcal{F}(\tau)^{sa}$. \\
\\
Before extending $\Phi$ to all of $\mathcal{F}(\tau)$, we show that $\Phi$ is positive and square-preserving.

\begin{lem} \label{L15.3 C1} $\Phi$ is positive and if $x\in \mathcal{F}(\tau)^{sa}$, then $\Phi(x^2)=\Phi(x)^2$.
\begin{proof} Since $\Phi$ is real linear and maps projections onto projections, it follows that $\Phi(x)\geq 0$, whenever $x\in \mathcal{G}_f^+$. If $x\in \mathcal{F}(\tau)^+$, then approximating $x$ by positive elements in $\mathcal{G}_f$ (see Proposition \ref{RP1 19/11/15}), using the continuity of $\Phi$ and closedness of  $\B^+$ (\cite[Theorem 4.2.2]{key-K1}), we obtain $\Phi(x)\geq 0$.

To show that $\Phi$ is square-preserving on $\mathcal{F}(\tau)^{sa}$, we start by noting that $\Phi$ is square-preserving on $\mathcal{G}^{sa}_f$, by Lemma \ref{nL1 04/12/15}. If $x\in \mathcal{F}(\tau)^{sa}$, we can once again use an approximation argument (and the joint continuity of multiplication in $\A$ and $\B$) to show that $\Phi(x^2)=\Phi(x)^2$. 
\end{proof}
\end{lem}

Finally, since any element $x\in \mathcal{F}(\tau)$ has a unique decomposition $x=x_1+ix_2$, where $x_1,x_2 \in \mathcal{F}(\tau)^{sa}$, we can extend $\Phi$ to a well-defined complex linear map on $\mathcal{F}(\tau)$ (with the properties claimed in Theorem \ref{nP extension2}). We have therefore proven the following result.

\begin{thm}\label{nP extension2}
Suppose $\Phi:\paf \ra \pb$ is a map such that $\Phi(p+q)=\Phi(p)+\Phi(q)$ whenever $p,q\in \paf$ with $pq=0$. If there exists a linear map $U$ from $\mathcal{F}(\tau)$ into $\xmx{\nu}{\B}$ such that $\Phi(p)=s(U(p))$ for all $p\in \paf$, and which has the property that $U(x_n)\rax{\mtop} U(x)$ whenever $\seq{x}\cup \{x\}\subseteq \mathcal{F}(\tau)$ is such that $x_n \ra x$ and $s(x_n)\leq s(x)$ for all $n\in \mathbb{N}^+$, then $\Phi$ can be extended to a positive linear map (still denoted $\Phi$) from $\mathcal{F}(\tau)$ into $\B$ such that $\norm{\Phi(x)}_\B\leq\norm{x}_\A$ and $\Phi(x^2)=\Phi(x)^2$ for all $x\in \mathcal{F}(\tau)^{sa}$.
\end{thm}

\section{Projection-preserving maps on $\mathcal{F}(\tau)$ that can be extended to Jordan $*$-homomorphisms} \label{extension}\label{S4}

In this section we show that the extension of the projection mapping obtained in the previous section can be  extended to a normal Jordan $*$-homomorphism from $\A$ into $\B$, provided the extension obtained in the previous section is normal.  We also wish to show that a normal projection-preserving linear mapping from $\mathcal{F}(\tau)$ into $\B$ can likewise be extended to a normal Jordan $*$-homomorphism from $\A$ into $\B$. It follows from the next result that these tasks can be accomplished concurrently.

\begin{thm} \label{nT3.6 L99}
If $\Phi: \mathcal{F}(\tau) \ra \B$ is a continuous linear operator, then the following are equivalent:
\begin{enumerate}
\item $\Phi$ maps projections in $\mathcal{F}(\tau)$ onto projections in $\B$;
\item $\Phi$ is adjoint-preserving and $\Phi(xy+yx)=\Phi(x)\Phi(y)+\Phi(y)\Phi(x)$ for every $x,y\in \mathcal{F}(\tau)$;
\item $\Phi(x)\geq 0$ and $\Phi(x^2)=\Phi(x)^2$ for any $x\in \mathcal{F}(\tau)^+$;
\item $\Phi(x)\geq 0$ for every $x\in \mathcal{F}(\tau)^+$, $\Phi(x^2)=\Phi(x)^2$ for any $x\in \mathcal{F}(\tau)^{sa}$ and $\norm{\Phi(x)}_\B \leq \norm{x}_\A$ for every $x\in \mathcal{F}(\tau)^{sa}$.
\end{enumerate}
\begin{proof}
The equivalence of the first three statements can be shown using a similar technique to the one employed in the proof of \cite[Theorem 3.6]{key-L99}. Furthermore, it is clear that (4) implies (3). To see that (3) implies (4), we note that if $x\in \mathcal{F}(\tau)^{sa}$, then $x=x_+-x_-$, with $x_+,x_-\in \mathcal{F}(\tau)^+$. Expanding $\Phi(x^2)$ using $x=x_+-x_-$, the linearity of $\Phi$ and the fact that $\Phi$ preserves Jordan products and products of positive elements yields $\Phi(x^2)=\Phi(x)^2$. It now follows from Proposition \ref{P18 16/10/14} that $\norm{\Phi(x)}_\B \leq \norm{x}_\A$ for every $x\in \mathcal{F}(\tau)^{sa}$.
\end{proof}
\end{thm}

Throughout this section we will let $\Phi: \mathcal{F}(\tau) \ra \B$ denote a linear, positive and normal map such that $\Phi(x^2)=\Phi(x)^2$ for every $x\in \mathcal{F}(\tau)^{sa}$. Note that it follows from Proposition \ref{P18 16/10/14}(8) that  $\Phi$ is automatically continuous. In order to extend $\Phi$, we start by defining a map $\Phi_1$ on projections. If $p\in\pa$, then it follows from the semi-finiteness of $\A$ that $\Dx{p}:=\{q\in \paf: q\leq p\}$ is an increasing net whose supremum is $p$. For the sake of convenience we will use an indexing set and let $\Dx{p}=\{p_\alpha:\alpha\in I_p\}$. Note that $\Phi$ is positive and therefore $\{\Phi(p_{\alpha})\}_{\alpha \in I_p}$ is increasing. By Proposition \ref{P18 16/10/14}(5), $\Phi(p_{\alpha})$ is a projection for every $\alpha\in I_p$ and so the SOT-limit of $\{\Phi(p_{\alpha})\}_{\alpha\in  I_p}$ exists and is a projection.  We define $\Phi_1(p)$ to be this limit. Note that $\Phi_1(p)=\Phi(p)$ for $p\in \paf$. 

In order to apply Lemma \ref{nL1 04/12/15} to extend $\Phi_1$ to $\mathcal{G}^{sa}$ we need to show that $\Phi_1$ is additive on orthogonal projections. Suppose $p$ and $q$ are orthogonal projections in $\pa$. Note that if $p_\alpha\in \Dx{p}$ and $q_\beta \in \Dx{q}$, then $p_\alpha + q_\beta \in \Dx{p+q}$, since $p_\alpha$ and $q_\beta$ are orthogonal projections, each with finite trace. Furthermore $p_\alpha+q_\beta \uparrow p+q$ and so if $r\in \Dx{p+q}$, then there exists an $(\alpha,\beta)\in I_p \times I_q$ such that $p_\alpha+q_\beta\geq r$.  It follows, using the positivity of $\Phi$,  that $\{\Phi(p_\alpha+q_\beta):\alpha \in I_p,\beta\in I_q\}$ is a subnet of $\{\Phi(r):r\in \Dx{p+q}\}$ and hence $\Phi(p_\alpha)+\Phi(q_\beta)=\Phi (p_\alpha+q_\beta)\sotc \Phi_1(p+q)$, since $\Phi_1(p+q):=\limm{r\in\Dx{p+q}}\Phi(r)$. However, we also have that $\Phi(p_\alpha)+\Phi(q_\beta)\sotc \Phi_1(p)+\Phi_1(q)$, by definition of $\Phi_1(p)$ and $\Phi_1(q)$ and the SOT-continuity of addition. Therefore $\Phi_1(p+q)=\Phi_1(p)+\Phi_1(q)$. By \cite[Exercise 2.3.4]{key-Con07} this further implies that $\Phi_1(p)$ and $\Phi_1(q)$ are orthogonal if $pq=0$. 

For $x=\summ{i=1}{n}\alpha_i p_i \in \mathcal{G}^{sa}$, define $\Phi_2(x):=\summ{i=1}{n}\alpha_i \Phi_1(p_i)$. It follows by Lemma \ref{nL1 04/12/15} that $\Phi_2$ is well-defined and real linear on commuting operators in $\mathcal{G}^{sa}$. Furthermore, $\Phi_2$ is a square-preserving contraction on $\mathcal{G}^{sa}$. Similarly to before, one can therefore use Proposition \ref{RP1 19/11/15} to extend $\Phi_2$ to a contraction $\Phi_3: \A^{sa}\ra \B$. Since $\Phi_1(p)=\Phi(p)$ for each $p\in \paf$, we have that $\Phi_2(x)=\Phi(x)$ for every $x\in \mathcal{G}_f^{sa}$ and therefore $\Phi_3$ agrees with $\Phi$ on $\mathcal{F}(\tau)^{sa}$. We will therefore also denote this extension $\Phi$. Note that $\Phi$ is real linear on commuting elements, but need not be linear on all elements of $\A^{sa}$. Recall that in section \ref{extension2} we showed that the existence of a linear map $U:\mathcal{F}(\tau)\ra \xmx{\nu}{\B}$ such that $\Phi(p)=s(U(p))$ for all $p\in \paf$ could be used, under certain circumstances, to show that the extension is linear. A similar idea is not applicable in this context, since the maps $U$ we will be interested in for applications of these results will typically not be defined on projections with infinite trace. We can however use the normality of $\Phi$ to prove that $\Phi(x)$ can be written as a SOT-limit for any $x\in \A^{sa}$ and hence demonstrate the real linearity of $\Phi$ on $\A^{sa}$. We need the following preliminary result.

\begin{lem} \label{L13.5a C3.1} \label{L13.5b C3.1}
If $p\in \pa$, then 
\begin{enumerate}
\item $\Phi(qpq)=\Phi(q)\Phi(p)\Phi(q)$ for all $q\in \paf$;
\item $\Phi(\id)\Phi(p)\Phi(\id)=\Phi(p)$.
\end{enumerate}
\begin{proof}
(1) Note that $p_\alpha \uparrowx{\alpha\in I_p} p$ implies that $qp_\alpha q \uparrowx{\alpha\in I_p} qpq$ and therefore $ \Phi(qp_\alpha q)\uparrowx{\alpha\in I_p} \Phi(qpq)$, using \cite[Proposition 1(vi)]{key-Dodds14} and the normality of $\Phi$. By \cite[Lemma 5.1.4]{key-K1}, this implies that $\Phi(qp_\alpha q)\sotcx{\alpha\in I_p} \Phi(qpq)$. Furthermore, by Proposition \ref{P18 16/10/14}(4), \cite[p.115]{key-K1} and the definition of $\Phi$, we obtain 
\begin{eqnarray*}
\limm{\alpha \in I_p}\Phi(qp_\alpha q)=\limm{\alpha\in I_p}\Phi(q)\Phi(p_\alpha)\Phi(q)=\Phi(q)(\limm{\alpha\in I_p}\Phi(p_\alpha)) \Phi(q)=\Phi(q)\Phi(p)\Phi(q).  
\end{eqnarray*} Combining these yields $\Phi(qpq)=\Phi(q)\Phi(p)\Phi(q)$, since SOT-limits are unique. \\

(2) If $\alpha \in I_p$, then $p_\alpha \in \Dx{\id}\equiv \D$ and so $\Phi(p_\alpha)\leq \supu{q\in \D}\Phi(q)=\Phi(\id)$. It follows that $\Phi(p_\alpha)=\Phi(\id)\Phi(p_\alpha)\Phi(\id) \sotcx{\alpha\in I_p} \Phi(\id)\Phi(p)\Phi(\id)$, by \cite[Proposition 2.5.2(3)]{key-Con07}, \cite[Proposition 1(vi)]{key-Dodds14} and \cite[Lemma 5.1.4]{key-K1}.
However, $\Phi(p_\alpha)\sotcx{\alpha\in I_p} \Phi(p)$ and so $\Phi(\id)\Phi(p)\Phi(\id)=\Phi(p)$.
\end{proof}
\end{lem}

We show that for any self-adjoint element $x$, $\Phi(x)$ can be viewed as a $SOT$-limit.

\begin{lem} \label{L13.6.0 C3.1} If $x\in \A^{sa}$, then $pxp\in \mathcal{F}(\tau)$ for every $p\in \paf$ and \[\Phi(x)=\limm{p\in \D} \Phi(pxp).\]
\begin{proof} Suppose $x=q\in \pa$. Then  
\begin{eqnarray}
\Phi(pqp)=\Phi(p)\Phi(q)\Phi(p) \sotcx{p\in \D}\Phi(\id)\Phi(q)\Phi(\id) =\Phi(q), \label{e13.5.1 C3.1}
\end{eqnarray}
using Lemma \ref{L13.5a C3.1} and the fact that $\Phi(p)\sotcx{p\in \D} \Phi(\id)$. If $x=\summ{i=1}{n}\alpha_i q_i\in \mathcal{G}^{sa}$, then for $p\in \paf$
\begin{eqnarray}
\Phi(pxp)=\summ{i=1}{n}\alpha_i\Phi(pq_i p) \sotcx{p\in \D}\summ{i=1}{n}\alpha_i\Phi(q_i) =\Phi(x), \label{e13.6a C3.1}
\end{eqnarray}
using the linearity of $\Phi$ on $\mathcal{F}(\tau)$, (\ref{e13.5.1 C3.1}) and the fact that the SOT is a vector topology. Next, suppose that $x\in \A^{sa}$. Let $\seq{x}\subseteq \mathcal{G}^{sa}$ be the sequence described in Proposition \ref{RP1 19/11/15} such that $x_n \ra x$ (and hence $\Phi(x_n)\ra \Phi(x)$). Let $\epsilon>0$ and fix $\eta \in K$, with $\norm{\eta}_K=1$. It follows that there exist  $n_0,n_1\in \n^+$, such that  
\begin{eqnarray}
\norm{x_n-x}_\A&<&\frac{\epsilon}{3} \label{e13.6.2 C3.1} \qquad \text{for all $n\geq n_0$ and} \\
\norm{\Phi(x_n)-\Phi(x)}_\B&<&\frac{\epsilon}{3} \qquad \text{for all $n\geq n_1$.} \label{e13.6.1 C3.1} 
\end{eqnarray}
Furthermore, for $p\in \paf$ and $n\geq \text{max}\{n_0,n_1\}$, we have 
\begin{eqnarray}
\norm{(\Phi(p(x-x_n)p))\eta}_K\leq \norm{p(x-x_n)p}_\A \norm{\eta}_K <\frac{\epsilon}{3}, \label{e13.6.3 C3.1} 
\end{eqnarray}
by Proposition \ref{P18 16/10/14}(8) and (\ref{e13.6.2 C3.1}). Fix $n\geq \text{max}\{n_0,n_1\}$. By (\ref{e13.6a C3.1}), we have that $\Phi(x_n)=\limm{p\in \D}\Phi(px_np)$. We can therefore find $q\in \paf$ such that $p\in \paf$ and $p\geq q$ implies that
\begin{eqnarray}
\norm{(\Phi(px_np)-\Phi(x_n))\eta}_K&<&\frac{\epsilon}{3}. \label{e13.6.4 C3.1}
\end{eqnarray}
Then for $p\in \paf$ with $p\geq q$, application of (\ref{e13.6.3 C3.1}), (\ref{e13.6.4 C3.1}) and (\ref{e13.6.1 C3.1}) yields
\begin{eqnarray*}
\norm{(\Phi(pxp)-\Phi(x))\eta}_K &\leq& \norm{(\Phi(pxp)-\Phi(px_np))\eta}_K + \\
&&\norm{(\Phi(px_np)-\Phi(x_n))\eta}_K + \norm{(\Phi(x_n)-\Phi((x))\eta}_K <\epsilon.
\end{eqnarray*}
It follows that $\Phi(pxp)\eta \raxx{p\in\D} \Phi(x)\eta$ for any $\eta\in K$ with $\norm{\eta}_K=1$ and therefore also for any $0\neq\eta\in K$ by considering $\tilde{\eta}=\frac{\eta}{\norm{\eta}}$. Therefore $\Phi(x)=\limm{p\in \D}\Phi(pxp)$.
\end{proof}
\end{lem}

Using Lemma \ref{L13.6.0 C3.1} (repeatedly), the fact that the SOT is a vector topology and the linearity of $\Phi$ on $\mathcal{F}(\tau)$, we can show that 
\[\alpha\Phi(x)+\beta \Phi(y)=\alpha \limm{p\in \D}\Phi(pxp)+ \beta\limm{p\in \D}\Phi(pyp) =\Phi(\alpha x +\beta y),\]  
whenever $x,y\in \A^{sa}$ and $\alpha, \beta \in \mathbb{R}$. We have therefore shown the real linearity of $\Phi$ on $\A^{sa}$. As in the previous section we can now extend $\Phi$ to a complex linear map on $\A$. Furthermore,  
\begin{eqnarray}
\Phi(x)=\limm{p\in \D}\Phi(pxp) \qquad \forall x\in \A. \label{eL13.6 C3.1}
\end{eqnarray}

\begin{lem} \label{L14.1 C3.1}\label{L14.3 C3.1}
$\Phi$ is a normal Jordan $*$-homomorphism from $\A$ into $\B$.
\begin{proof}
Recalling the definition of $\Phi_1$ and noting that $\Phi$ is an extension of $\Phi_1$, it follows that $\Phi$ maps projections onto projections. Furthermore, if $x\in \A$, then $x=x_1+ix_2$ for some $x_1,x_2\in \A^{sa}$ and so  
\begin{eqnarray*}
\norm{\Phi(x)}_\B\leq\norm{\Phi(x_1)}_\B+\norm{\Phi(x_2)}_\B \leq \norm{x_1}_\A+\norm{x_2}_\A \leq 2\norm{x}_\A,
\end{eqnarray*}
since $\Phi$ is a contraction on self-adjoint elements. It follows that $\Phi$ is continuous (and projection-preserving) and hence a Jordan $*$-homomorphism, by \cite[Theorem 3.7]{key-L99}. Next, we show that $\Phi$ is normal. Suppose $\{x_\gamma\}_\gamma \cup \{x\}\subseteq \A^+$ such that $x_\gamma \uparrow x$. In this case, we have that for any $p\in \paf$, $px_\gamma p\uparrow pxp$, by \cite[Proposition 1(vi)]{key-Dodds14}. Therefore $\Phi(px_\gamma p)\uparrow \Phi(pxp)$, since $\Phi$ is normal on $\mathcal{F}(\tau)$. It follows by \cite[Lemma 5.1.4]{key-K1} that
\begin{eqnarray}
\Phi(pxp)&=&\limm{\gamma}\Phi(px_\gamma p). \label{e14.1.1 C3.1}
\end{eqnarray}
Since $\Phi$ is positive, we have that $\{\Phi(x_\gamma)\}_\gamma$ is increasing and $0\leq\Phi(x_\gamma)\leq \Phi(x) \leq \norm{\Phi(x)}_\B \id$ for all $\gamma$, by \cite[Proposition 4.2.3]{key-K1}. Therefore $\limm{\gamma}\Phi(x_\gamma)$ exists, by \cite[Lemma 5.1.4]{key-K1}. Furthermore, for any $p\in \paf$,
\[\Phi(p)\Phi(x)\Phi(p)=\Phi(pxp) =\limm{\gamma}\Phi(px_\gamma p) =\limm{\gamma}[\Phi(p)\Phi(x_\gamma)\Phi(p)] =\Phi(p)[\limm{\gamma}\Phi(x_\gamma)]\Phi(p),\] by \cite[Exercise 10.5.21(ii)]{key-K4} and (\ref{e14.1.1 C3.1}). We also have that $\Phi(p)\sotcx{p\in \D} \Phi(\id)$ and so $\Phi(p)\Phi(x)\Phi(p)\sotcx{p\in \D}\Phi(\id)\Phi(x)\Phi(\id) =\Phi(x)$, by \cite[Exercise 10.5.21(ii)]{key-K4}. Similarly, $\Phi(p)[\limm{\gamma}\Phi(x_\gamma)]\Phi(p)\sotcx{p\in \D}\limm{\gamma}\Phi(x_\gamma)$. By combining these facts, we obtain $\Phi(x)= \limm{\gamma}\Phi(x_\gamma)$. Since, $\limm{\gamma}\Phi(x_\gamma)=\supu{\gamma}\Phi(x_\gamma)$, by \cite[Lemma 5.1.4]{key-K1}, we therefore have that $\Phi(x_\gamma)\uparrow \Phi(x)$.
\end{proof}
\end{lem}

We conclude by showing that the extension we obtained in this section is unique. For this exercise, let $\Phi_0$ denote the original map from $\mathcal{F}(\tau)$ into $\B$ and suppose $\Psi:\A \ra \B$ is another normal Jordan $*$-homomorphism extending $\Phi_0$. Let $x\in \A$. Since $\Psi$ is normal, we have that $\Psi(p)\uparrowx{p\in \D} \Psi(\id)$ and hence it follows from \cite[Lemma 5.1.4]{key-K1} that $\Psi(p) \sotcx{p\in \D} \Psi(\id)$. Therefore, $\Psi(pxp)= \Psi(p)\Psi(x)\Psi(p) \sotcx{p\in \D} \Psi(\id)\Psi(x)\Psi(\id) =\Psi(x)$, by \cite[Exercise 10.5.21]{key-K4}. Similarly, $\Phi(pxp) \sotcx{p\in \D} \Phi(x)$, but $\Phi(pxp)=\Phi_0(pxp)=\Psi(pxp)$ for every $p\in \paf$ and so $\Phi(x)=\Psi(x)$. Since this holds for every $x\in \A$, we have that $\Psi =\Phi$. We summarize the results of this section in the following theorem.

\begin{thm} \label{P extension}
Suppose $(\A,\tau)$ and $(\B,\nu)$ are semi-finite von Neumann algebras. If $\Phi: \mathcal{F}(\tau) \ra \B$ is linear, positive, normal and square-preserving on self-adjoint elements (equivalently, continuous, linear, normal and projection-preserving), then $\Phi$ can be extended uniquely to a normal Jordan $*$-homomorphism (also denoted $\Phi$) from $\A$ into $\B$. Furthermore, in this case, $\Phi(x)=\limx{p\in \D}\Phi(pxp)$ for all $x \in \A$, where the limit is taken in the SOT, and $\norm{\Phi(x)}_\B\leq\norm{x}_\A$ for all $x\in \A^{sa}$. 
\end{thm}

\begin{rem}\label{R1 03/07/18}
If, in addition to the assumptions of Theorem \ref{P extension}, $p=0$, whenever $p\in \paf$ with $\Phi(p)=0$, then the extension obtained by Theorem \ref{P extension} will be isometric on $\A^{sa}$. The reason is as follows. Recall that $\Phi$ is used to define $\Phi_1$ on projections.  If $p>0$, then there exists a $\beta\in I_p$ such that $0\neq p_\beta \in \Dx{p}$, since $\A$ is semi-finite. $\{\Phi(p_\alpha)\}_{\alpha \in I_p}$ is an increasing net and so $\Phi_1(p)\geq \Phi(p_\beta)>0$. It follows by Remark \ref{R1 01/07/18} that $\Phi_2$, the extension of $\Phi_1$ to $\mathcal{G}$, is isometric on self-adjoint elements. Since $\mathcal{G}^{sa}$ is dense in $\A^{sa}$ it follows that the extension of $\Phi_2$ to $\A^{sa}$ is isometric.
\end{rem}

\section{Projection mappings that can be extended to Jordan $*$-homomorphisms}\label{S5}

Thus far we have extended a projection mapping firstly to a linear map from $\mathcal{F}(\tau)$ into $\B$ (using the linearity of another map $U$) and then subsequently to a Jordan $*$-homomorphism from $\A$ into $\B$, provided the first extension is normal on $\mathcal{F}(\tau)$. Our motivation was two-fold. Firstly, there are occasions (the characterization of non-commutative composition operators, for example) where one is interested in extending a projection-preserving mapping which is  defined on all of $\mathcal{F}(\tau)$ and not just on projections with finite trace. The second reason is that there are situations where the normality of the first extension can be shown using particular properties of the map $U$ or the space on which $U$ is defined. In this section, however, we provide sufficient conditions for a projection mapping, defined on the lattice of projections with finite trace, to be extended directly to a Jordan $*$-homomorphism on the whole von Neumann algebra.\\

\begin{thm}\label{T2 04/07/18}\label{L2 03/07/18}
Suppose $\Phi:\paf \ra \pb$ is a map such that  $\Phi(p+q)=\Phi(p)+\Phi(q)$ whenever $p,q\in \paf$ with $pq=0$. Suppose further there exists a positive normal linear map $U$ from $\mathcal{F}(\tau)$ into $\xmx{\nu}{\B}$ such that $\Phi(p)=s(U(p))$ for all $p\in \paf$, and which has the properties
\begin{itemize}
\item $U(x_n)\rax{\mtop} U(x)$ whenever $\seq{x}\cup \{x\}\subseteq \mathcal{F}(\tau)$ is such that $x_n \ra x$ and $s(x_n)\leq s(x)$ for all $n\in \mathbb{N}^+$,
\item $\nu(s(U(p)))<\infty$ for every $p\in \paf$ (equivalently $\Phi(\paf)\subseteq \pbf$).
\end{itemize}
Let $\Phi:\mathcal{F}(\tau)\ra\B$ denote the extension obtained in Theorem \ref{nP extension2}. 
\begin{enumerate}
\item If $x\in \mathcal{F}(\tau)^{sa}$ and $p\in \paf$ with $p\geq s(x)$, then  $\Phi(x)U(p)=U(x)=U(p)\Phi(x)=U(p)^{1/2}\Phi(x)U(p)^{1/2}$;
\item If $x\in \mathcal{F}(\tau)^{sa}$ and $p\in \paf$ with $p\geq s(x)$, then there exists a $v_p\in \nmb$ such that $U(p)^{1/2}v_p=\Phi(p)=v_pU(p)^{1/2}$ and $\Phi(x)=v_pU(x)v_p$;
\item  $\Phi$ can be extended to a normal Jordan $*$-homomorphism (still denoted $\Phi$) from $\A$ into $\B$. Furthermore, in this case, $\Phi(x)$ is the SOT-limit  of $\{\Phi(pxp)\}_{p\in \D}$ for any $x \in \A$ and $\norm{\Phi(x)}_\B\leq \norm{x}_\A$ for all $x\in \A^{sa}$.
\end{enumerate}
\begin{proof} 
(1) The second equality is given by Lemma \ref{L15.1 C1}. Since, in addition, $U$ is positive, we have that $r(U(p))=s(U(p))=\Phi(p)$, in this case. A similar argument to the one employed to prove Lemma \ref{L15.1 C1} can therefore be used to obtain the first equality. Since $\Phi(x)$ and $U(p)$ commute and $\Phi(x)\in \B$, the same is true of $\Phi(x)$ and $U(p)^{1/2}$, by \cite[Theorem 5.12]{key-dP} and \cite[Corollary 15]{key-Terp1}. Using this commutativity, we obtain the third equality.

(2) We start by showing the existence of $v_p\in \nmb$ such that $U(p)^{1/2}v_p=\Phi(p)=v_pU(p)^{1/2}$. Let $\M$ denote the reduced von Neumann algebra generated by $\Phi(p)\B\Phi(p)$ and let $\varphi$ denote the isometric $*$-isomorphism from $\Phi(p)\B\Phi(p)$ onto $\M$. (See \cite[p.211,212]{key-Dodds14} for details regarding reduced spaces.) Note that $\varphi$ extends to a $*$-isomorphism from $\Phi(p)S(\B,\nu)\Phi(p)$ onto $S(\M)=S(\M,\tau_{\M})$, the set of all trace-measurable operators affiliated with $\M$, equipped with the reduced trace $\tau_\M$. Since $\nu(s(U(p)))<\infty$, $\M$ is a trace-finite von Neumann algebra. Furthermore, since $s(U(p)^{1/2})=s(U(p))=\Phi(p)$ and $\varphi(\Phi(p))$ is the identity of $\M$, $\varphi(U(p)^{1/2})$ has a closed densely defined inverse. It is easily checked that this inverse is affiliated with $\M$. It follows that $\varphi(U(p)^{1/2})$ is invertible in $S(\M)$ using the trace-finiteness of $\M$. Let $v$ denote the inverse of $\varphi(U(p)^{1/2})$. By \cite[Proposition 1]{key-Dodds14}, $v\geq 0$. It follows that $v_p:=\varphi^{-1}(v)\geq 0$ and  $U(p)^{1/2}v_p=\Phi(p)=v_pU(p)^{1/2}$.  Note that by Lemma \ref{L15.2 C1} $s(\Phi(x))=r(\Phi(x))\leq s(U(p))=\Phi(p)$. We therefore have that 
\[v_pU(x)v_p=v_p\left(U(p)^{1/2}\Phi(x)U(p)^{1/2}\right)v_p=\Phi(p)\Phi(x)\Phi(p)=\Phi(x).\]

(3) It was shown in  section \ref{extension2} that the map $\Phi:\mathcal{F}(\tau)\ra\B$ is a positive linear map such that $\norm{\Phi(x)}_\B\leq\norm{x}_\A$ and $\Phi(x^2)=\Phi(x)^2$ for all $x\in \mathcal{F}(\tau)^{sa}$. We show that $\Phi$ is normal. Suppose $\net{x}\cup \{x\}\subseteq \mathcal{F}(\tau)^+$ is such that $x_\lambda \uparrow x$. Let $p=s(x)$. Then $p\in \paf$, since $x\in \mathcal{F}(\tau)$. Furthermore, $p\geq s(x_\lambda)$ for each $\lambda$, since $x_\lambda \uparrow x$. Since $U$ is normal, we have that $U(x_\lambda) \uparrow U(x)$. Using \cite[Proposition 1]{key-Dodds14} and (2)  this implies that 
\[\Phi(x_\lambda)=v_pU(x_\lambda)v_p\uparrow v_pU(x)v_p=\Phi(x).\]
Since $\Phi$ is normal, it follows by Theorem \ref{P extension} that $\Phi$ can be extended to a normal Jordan $*$-homomorphism from $\A$ into $\B$, with the desired properties.
\end{proof}
\end{thm}

\section{Surjectivity}

We have seen that the result of the extension procedures we have detailed is a normal Jordan $*$-homomorphism. In this section we describe sufficient conditions to ensure that such a map is surjective. We need a preliminary result. To prove this result we will need some information regarding the ultra-weak operator topology (UWOT). The ultra-weak operator topology is the locally convex Hausdorff topology on $\bh$ generated by the semi-norms $\rho_{(\xi),(\eta)}$, where \[\rho_{(\xi),(\eta)}(x):=|\summ{n=1}{\infty} \ip{x\xi_n}{\eta_n}|\] and $(\xi)=\seq{\xi}$ and $(\eta)=\seq{\eta}$ are sequences in $H$ satisfying $\summ{n=1}{\infty}\norm{\xi_n}^2<\infty$ and $\summ{n=1}{\infty}\norm{\eta_n}^2<\infty$. The ultra-weak operator topology is stronger than the weak operator topology and the two topologies coincide on norm bounded sets. A positive linear mapping between von Neumann algebras is normal if and only if it is continuous with respect to the ultra-weak operator topologies (\cite{key-dP}). 

\begin{lem} \label{P9 16/10/14}
Suppose $\Phi:\A\ra \B$ is a normal Jordan $*$-homomorphism. If there exists a $k>0$ such that $\norm{x}\leq k\norm{\Phi(x)}$ for all $x\in \A$, then $\Phi(A)$ is WOT-closed.
\begin{proof} 
Since $\Phi$ is positive, linear and normal, it is therefore continuous with respect to the ultra-weak operator topologies on $\A$ and $\B$. Suppose $b\in \overline{\Phi(\A)}^{WOT}$ with $\norm{b}=1$. It follows from the Kaplansky density theorem (\cite[Theorem 5.3.5]{key-K1}) that there exists a net $\net{b}\subseteq \Phi(\A)$ such that $\norm{b_\lambda}\leq 1$ for each $\lambda$ and $b_\lambda \rax{WOT} b$. Let $\net{a}\subseteq \A$ be such that $\Phi(a_\lambda)=b_\lambda$ for each $\lambda$. Then \[\norm{a_\lambda}\leq k \norm{\Phi(a_\lambda)}=k \norm{b_\lambda} \leq k\] for each $\lambda$. The unit ball (and hence any multiple of the unit ball) in $\A$ is WOT-compact and so there exists a subnet $\netp{a}\subseteq \net{a}$ such that $a_{\lambda'}\rax{WOT} a$ for some $a\in \A$ with $\norm{a}\leq k$. Since the weak operator topology and ultra-weak operator topology coincide on norm-bounded sets, $a_{\lambda'}\rax{UWOT} a$. Using the continuity of $\Phi$, we have that $\Phi(a_{\lambda'})\rax{UWOT} \Phi(a)$. It follows that $\Phi(a_{\lambda'})\rax{WOT} \Phi(a)$, since the ultra-weak operator topology is stronger than the weak operator topology. However, $\Phi(a_{\lambda'})\rax{WOT} b$ and so $\Phi(a)=b$, since the weak operator topology is Hausdorff.
\end{proof}
\end{lem}

We are now in a position to describe conditions which ensure that a normal Jordan $*$-homomorphism is surjective.

\begin{prop} \label{P19 16/10/14}
Suppose $\Phi:\A\ra \B$ is a unital normal Jordan $*$-homomorphism which is isometric on $\A^{sa}$. If $\Phi(p)\B \Phi(p)\subseteq \Phi(\A)$ for all $p\in \paf$, then $\Phi$ is a Jordan $*$-isomorphism.
\begin{proof}
Suppose $b\in \B$. Since $\Phi$ is normal and unital, we have that $\Phi(p)\uparrow_{p\in \D} \Phi(\id)=\id$. It therefore follows from \cite[Lemma 5.1.4]{key-K1} that $\Phi(p)\sotc \id$. Note that $\Phi(p)b\Phi(p) \sotc \id b \id =b$ and therefore $\Phi(p)b\Phi(p) \rax{WOT} b$. Since $\Phi(p)b\Phi(p)\in \Phi(\A)$ for all $p\in \paf$, it follows that  $b\in\overline{\Phi(\A)}^{\text{WOT}}$. 

We wish to apply Lemma \ref{P9 16/10/14} to show that $\overline{\Phi(\A)}^{\text{WOT}}=\Phi(\A)$ and hence that $b\in \Phi(\A)$. We therefore show that there exists a $k>0$ such that $\norm{x}\leq k\norm{\Phi(x)}$ for all $x\in \A$. Let $x\in \A$. Then $x=x_1+ix_2$ for some $x_1,x_2\in \A^{sa}$. $\Phi$ is Jordan $*$-homomorphism and so $\Phi(x^*)=\Phi(x)^*$ for all $x\in \A$. It follows that $\Phi(x_1)$ and $\Phi(x_2)$ are the real and imaginary parts of $\Phi(x)$. Furthermore,
\begin{eqnarray*}
\norm{x}_\A&=&\norm{x_1+ix_2}_\A \leq\norm{x_1}_\A+\norm{x_2}_\A =\norm{\Phi(x_1)}_\B+\norm{\Phi(x_2)}_\B  \\
&=&\norm{\text{Re}(\Phi(x))}_\B+\norm{\text{Im}(\Phi(x))}_\B \leq\norm{\Phi(x)}_\B+\norm{\Phi(x)}_\B=2 \norm{\Phi(x)}_\B, 
\end{eqnarray*}
where we have used the fact that $\Phi$ is isometric on $\A^{sa}$. This shows that $\Phi$ is injective and allows application of Lemma \ref{P9 16/10/14} to show that $\overline{\Phi(\A)}^{\text{WOT}}=\Phi(\A)$ and hence that $b\in \Phi(\A)$.
\end{proof}
\end{prop}

\section*{Conclusion}

The conditions required for the extension procedures in this article may seem restrictive, but it will be shown that they are sufficiently general to facilitate structural descriptions of isometries between various non-commutative spaces (\cite{key-dJ1},\cite{key-dJ2} and \cite{key-dJ3}) and are useful in the characterization of quantum composition operators between symmetric spaces associated with semi-finite von Neumann algebras (\cite{key-dJ18c}).

\section*{Acknowledgements}

The majority of this research was conducted during the first author's doctoral studies at the University of Cape Town. The first author would like to thank his Ph.D. supervisor, Dr Robert Martin, for his input and guidance.

\bibliographystyle{amsplain}

\begin{thebibliography}{00}
\bibitem{key-Bunce92} L.J. Bunce and J.D.M. Wright, \emph{The Mackey-Gleason problem}, Bull. Amer. Math. Soc. (N.S.)  \textbf{26} (1992), 288-293.
\bibitem{key-Bunce93} L.J. Bunce and J.D.M. Wright, \emph{On Dye's theorem for Jordan operator algebras}, Expo. Math. \textbf{11} (1993), 91-95.
\bibitem{key-Bunce94} L.J. Bunce and J.D.M. Wright, \emph{The Mackey-Gleason problem for vector measures on projections in von Neumann algebras}, J. London. Math. Soc. \textbf{49}(2) (1994), 133-149.
\bibitem{key-Con07}J.B. Conway, \emph{A course in functional analysis, Second edition}, Springer, 2007.
\bibitem{key-dJ17} P. de Jager, \emph{Isometries on symmetric spaces associated with semi-finite von Neumann algebras}, Thesis (2017), (Available online at \url{https://open.uct.ac.za/handle/11427/25167})
\bibitem{key-dJ18c}P. de Jager, \emph{Quantum composition operators on symmetric spaces associated with semi-finite von Neumann algebras} (in preparation).
\bibitem{key-dJ1} P. de Jager and J.J. Conradie, \emph{Isometries on quantum symmetric spaces associated with semi-finite von Neumann algebras}(in preparation).
\bibitem{key-dJ2} P. de Jager and J.J. Conradie, \emph{Isometries on quantum Lorentz spaces associated with semi-finite von Neumann algebras} (in preparation).
\bibitem{key-dJ3} P. de Jager and J.J. Conradie, \emph{Isometries on quantum Orlicz spaces associated with semi-finite von Neumann algebras}, (in preparation).
\bibitem{key-Dodds14}P.G. Dodds and B. de Pagter, \emph{Normed K\"othe spaces: A non-commutative viewpoint}, Indag. Math. \textbf{25} (2014), 206-249 .
\bibitem{key-dP}P.G. Dodds, B. de Pagter and F.A. Sukochev, \emph{Theory of noncommutative integration}, unpublished monograph (to appear).
\bibitem{key-Dye55}H.A. Dye, \emph{On the geometry of projections in certain operator algebras}, Ann. of Math., Second Series, \textbf{61}(no.1) (1955), 73-89.
\bibitem{key-Ham15} J. Hamhalter, \emph{Dye's Theorem and Gleason's Theorem for $AW^*$-algebras}, J. Math. Anal. Appl. \textbf{422} (2015), 1103-1115.
\bibitem{key-K1}R.V. Kadison and J.R. Ringrose, \emph{Fundamentals of the theory of operator algebras, Volume 1}, Birkh\"auser, Academic Press, (1983).
\bibitem{key-K3}R.V. Kaddison and J.R. Ringrose, \emph{Fundamentals of the theory of operator algebras, Volume 3}, Birkh\"auser, Academic Press, (1983).
\bibitem{key-K4}R.V. Kaddison and J.R. Ringrose, \emph{Fundamentals of the theory of operator algebras, Volume 4: Advanced theory - An exercise approach}, Birkh\"auser, Academic Press, (1983).
\bibitem{key-L99}L.E. Labuschagne, \emph{Composition operators on non-commutative $L^p$-spaces}, Expo. Math.  \textbf{17} (1999), 429-468.
\bibitem{key-Terp1}M. Terp, \emph{$L^p$-spaces associated with von Neumann algebras} Notes, Copenhagen University, (1981).
\bibitem{key-Y81}F.J. Yeadon, \emph{Isometries of non-commutative $L^p$-spaces}, Math. Proc. Camb. Phil. Soc. \textbf{90} (1981), 41-50.
\end{thebibliography}

\end{document}